\documentclass[11pt, reqno]{amsart}
\usepackage{version}
\newtheorem{theorem}{Theorem}[section]

\newtheorem{prop}[theorem]{Proposition}

\theoremstyle{definition}
\newtheorem{definition}[theorem]{Definition}

\newtheorem{remark}[theorem]{Remark}

\newcommand{\Hol}{{\rm Hol}}

\newcommand{\B}{\mathbb{B}}
\newcommand{\C}{\mathbb{C}}

\newcommand{\N}{\mathbb{N}}

\newcommand{\R}{\mathbb{R}}

\newcommand{\Aut}{\mathop{{\rm Aut}}}

\def\v{\varphi}
\def\de{\partial}

\numberwithin{equation}{section}

\begin{document}
\title[Solving the Loewner PDE]{Solving the Loewner PDE in complete hyperbolic starlike domains of $\C^N$}
\author[L. Arosio]{Leandro Arosio$^\dag$}
\author[F. Bracci]{Filippo Bracci$^\dag$}
\author[E. F. Wold]{Erlend Forn\ae ss Wold$^{\dag\dag}$}
\address{F. Bracci, L. Arosio: Dipartimento Di Matematica\\
Universit\`{a} di Roma \textquotedblleft Tor Vergata\textquotedblright\ \\
Via Della Ricerca Scientifica 1, 00133 \\
Roma, Italy} \email{arosio@mat.uniroma2.it, fbracci@mat.uniroma2.it}
\address{E. F. Wold: Department of Mathematics\\
University of Oslo\\
Postboks 1053 Blindern, NO-0316 Oslo, Norway}\email{erlendfw@math.uio.no}

%
%    General info
%
\subjclass[2010]{}
\date{\today}
\thanks{$\dag$ Partially supported by the ERC grant ``HEVO - Holomorphic Evolution Equations'' n. 277691}
\thanks{$^{\dag\dag}$ Partially supported by the NFR grant 209751/F20}
\begin{abstract}
We prove that any Loewner PDE in a complete hyperbolic starlike domain of $\C^N$ (in particular in bounded convex domains)  admits an essentially unique univalent solution with values in $\C^N$.
\end{abstract}

\maketitle

\section{Introduction}

The classical Loewner theory in the unit disc of $\C$ was introduced by Ch. Loewner \cite{Loewner} in 1923 and later developed by P.P. Kufarev \cite{Kuf1943} and Ch. Pommerenke \cite{Pommerenke-65}. We refer the reader to \cite{ABCD} for a recent survey on  applications and generalizations of such a theory.

In higher dimension, J. Pfaltzgraff \cite{Pf74,Pf75} extended the basic theory to $\C^N$, and later on many authors contributed to study the higher (or even infinite) dimensional Loewner ODE and PDE. Just to name a few, we recall here the contributions of T. Poreda \cite{Por91}, I. Graham, H. Hamada and G. Kohr \cite{GHK1,GHK2}.

More recently, the second named author together with M. D. Contreras and S. D\'iaz-Madrigal \cite{BCD1}, \cite{BCD2} generalized and solved the Loewner ODE on complete hyperbolic manifolds, and later the first two authors together with H. Hamada and G. Kohr \cite{ABHK} showed that the Loewner PDE on complete hyperbolic manifolds always admits an (essentially unique) abstract univalent solution with values in a complex manifold. The main remaining
open problem, whether  any  Loewner PDE in a complete hyperbolic domain in $\C^N$ admits a univalent solution with values in $\C^N$, has been given many partial positive answers, see \cite{DGHK,GHK1,GHK2,Ar1,Ar2,Ar3,Voda}.

In the present paper we show that any Loewner PDE in a complete hyperbolic starlike domain in $\C^N$ (in particular in bounded convex domains and in the unit ball)  admits a univalent solution with values in $\C^N$. For $N=1$ it is known \cite{one} that any Loewner PDE in the unit disc admits a univalent solution with values in $\C$, so in what follows we will focus on  the case $N\geq 2$.

Referring the reader to Section \ref{chains} for definitions and comments, our main result can be stated in the following way:

\begin{theorem}\label{main-intro}
Let $D\subset \C^N$ be a complete hyperbolic starlike domain.
Let $G:D\times \R^+\to \C^N$ be a Herglotz vector
field of order $d\in[1,+\infty]$. Then there exists a family  of univalent mappings $(f_t\colon D\to \C^N)$ of order $d$ which solves the Loewner PDE
\begin{equation}\label{L-PDE}
\frac{\partial f_t}{\partial t}(z)=-df_t(z)G(z,t),\quad\mbox{a.a. }t\geq 0,\forall z\in D.
\end{equation}
Moreover, $R:=\cup_{t\geq 0}f_t(D)$ is a Runge and Stein domain in $\C^N$ and any other solution to \eqref{L-PDE} is of the form $(\Phi \circ f_t)$ for a suitable holomorphic map $\Phi:R\to \C^N$.
\end{theorem}

\section{Generalized Loewner theory}\label{chains}

Let $D\subset \C^N$ be a domain. Recall that a holomorphic vector field $H$ on $D$ is {\sl semicomplete} if the Cauchy problem $\stackrel{\cdot}{x}(t)=H(x(t))$, $x(0)=z_0$ has a solution defined for all $t\in [0,+\infty)$ for all $z_0\in D$. Semicomplete holomorphic vector fields on complete hyperbolic manifolds have been characterized in terms of the Kobayashi distance (see, {\sl e.g.}, \cite{AB}).

\begin{definition}\label{Herglotz}
Let $D\subset\C^N$ be a domain. A \textit{Herglotz vector field of order $d\geq 1$} on $D$ is a mapping $G:D\times \R^+\to
\C^N$ with the following properties:
\begin{itemize}
\item[(i)] The mapping $G(z,\cdot)$ is measurable on $\R^+$ for all
$z\in D$.
\item[(ii)] The mapping $G(\cdot,t)$ is a holomorphic vector field on $D$ for all $t\in \R^+$.
\item[(iii)] For any compact set $K\subset D$ and all $T>0$, there exists
a function $C_{T,K}\in L^d([0,T],\mathbb{R}^+)$ such that
$$\|G(z,t)\|\leq C_{T,K}(t),\quad z\in K, \mbox{ a.a.}\ t\in [0,T].$$
\item[(iv)]$D\ni z\mapsto G(z,t)$ is semicomplete
for almost all $t\in [0,+\infty)$.
\end{itemize}
\end{definition}

Herglotz vector fields are strictly related to evolution families:

\begin{definition}\label{EF}
Let $D\subset \C^N$ be a domain.
A family $(\v_{s,t})_{0\leq s\leq t}$ of holomorphic
self-mappings of $D$ is  an {\sl
 evolution family of order $d\geq [1,+\infty]$} if it satisfies
\begin{equation}\label{evolution_property}
\v_{s,s}={\sf id},\quad \v_{s,t}=\v_{u,t}\circ \v_{s,u},\quad 0\leq s\leq u\leq t,
\end{equation}
 and if for any $T>0$ and for any compact set
  $K\subset\subset D$ there exists a  function $c_{T,K}\in
  L^d([0,T],\R^+)$ such that
  \begin{equation}\label{ck-evd}
\|\v_{s,t}(z)-\v_{s,u}(z)\|\leq \int_{u}^t
c_{T,K}(\xi)d \xi, \quad z\in K,\  0\leq s\leq u\leq t\leq T.
  \end{equation}
\end{definition}
\begin{remark}
For all $0\leq s\leq t$ the mapping $\v_{s,t}\colon D\to D$ is univalent (see {\sl e.g.} \cite[Proposition 2.3]{ABHK})
\end{remark}

As a consequence of the main results in \cite{BCD2} and \cite{AB} we have the following solution to the generalized Loewner ODE:

\begin{theorem}\label{prel-thm}
Let $D$ be a complete hyperbolic domain in $\C^N$. Then for any  Herglotz vector field $G(z,t)$ of order
$d\in [1,+\infty]$ there exists a unique evolution family $(\v_{s,t})$ of order $d$
  over $D$ such that for all
$z\in D$
  \begin{equation}\label{solve}
\frac{\de \v_{s,t}}{\de t}(z)=G(\v_{s,t}(z),t) \quad
\hbox{a.a.\ } t\in [s,+\infty).
  \end{equation}
Conversely for any  evolution family $(\v_{s,t})$ of order $d\in [1,+\infty]$ over $D$ there exists a  Herglotz vector field $G$ of
order $d$ such that \eqref{solve} is satisfied. Moreover,
if $H$ is another weak holomorphic vector field which satisfies
\eqref{solve} then $G(z,t)=H(z,t)$ for all $z\in D$ and almost
every $t\in \R^+$.
\end{theorem}

We now define the Loewner PDE and its solutions.
\begin{definition}\label{solution}
Let $D$ be a complete hyperbolic domain in $\mathbb{C}^N$. The partial differential equation
\begin{equation}\label{LoewnerPDE}
\frac{\partial f_t(z)}{\partial t}=-df_t(z)G(z,t),
\end{equation}
 where $G(z,t)$ is a  Herglotz vector field of order $d\in[1,+\infty]$, is called the {\sl Loewner PDE}.
A {\sl solution} to (\ref{LoewnerPDE}) is a family of holomorphic mappings $(f_t)$ from $D$ to a complex manifold $Q$ of dimension $N$  such that
\begin{itemize}
\item[(i)]   the mapping  $t\mapsto f_t$ is continuous with respect to the topology  in $\Hol(D,Q)$ induced by the uniform convergence on compacta in $D$,
\item[(ii)] for all fixed $z\in D$ the mapping $t\mapsto f_t(z)$ is locally absolutely continuous in $\R^+$,
\item[(iii)] for all fixed $z\in D$ equality  (\ref{LoewnerPDE})  holds a.e. in  $\R^+$.
\end{itemize}
\end{definition}

The following proposition has been proved for $D=\mathbb{D}\subset \C$ in \cite[Proposition 2.3]{duality}, and the proof can be adapted  to
several variables (see also \cite{ABHK}).
\begin{prop}\label{urk}
Let $D\subset \C^N$ be a complete hyperbolic domain, let $Q$ be a complex manifold of dimension $N$ endowed with a Hermitian distance $d_Q$. Let $G(z,t)$ be a  Herglotz vector field on $D$ of order $d\in[1,+\infty]$.
Let  $(f_t\colon D\to Q)$ be a solution  to  the  Loewner PDE $(\ref{LoewnerPDE})$. Then
\begin{itemize}
\item[(i)] the family $(f_t\colon D\to\ Q)$ is of order $d$, that is for any compact set $K\subset D$ and all $T>0$, there exists
a function $c_{K,T}\in L^d([0,T],\mathbb{R}^+)$ such that
$$d_Q(f_s(z),f_t(z))\leq \int_{s}^t
c_{T,K}(\xi)d \xi, \quad z\in K,\  0\leq s\leq t\leq T,$$
\item[(ii)] there exists a set of zero measure $E\subset \R^+$ such that for all $z\in D$ and all $t\in \R^+\setminus E$ the partial  derivative $\partial f_t(z)/\partial t$ exists and equality (\ref{LoewnerPDE}) holds,
\item[(iii)]  if $(\v_{s,t}\colon D\to D)$ is the  evolution family which solves the Loewner ODE (\ref{solve}), then the following functional equation holds:
\begin{equation}\label{functional}
 f_s=f_t\circ \v_{s,t},\quad 0\leq s\leq t.
 \end{equation}
\end{itemize}
\end{prop}
\begin{remark}\label{cresce}
Any family of univalent mappings  $(f_t\colon D\to Q)$ satisfying the functional equation (\ref{functional})
has growing images:
\begin{equation}\label{cch}
f_s(D)\subseteq f_t(D),\quad 0\leq s\leq t.
\end{equation}
\end{remark}

In fact, a family $(f_t\colon D\to Q)$ of univalent mappings satisfying \eqref{cch} and i) of Proposition \ref{urk} is called a {\sl Loewner chain of order $d$}. By the results in \cite{ABHK}, given any Loewner chain  $(f_t\colon D\to Q)$  of order $d\in[1,+\infty]$ there exists an Herglotz vector field $G(z,t)$ of order $d$ on $D$ such that  $(f_t)$ solves the Loewner PDE  (\ref{LoewnerPDE}).

The following theorem follows from various results in  \cite{ABHK} :
\begin{theorem}\label{nostro}
Let $D\subset \C^N$ be  a complete hyperbolic domain and let $G(z,t)$ be  a  Herglotz vector field of order $d\in[1,+\infty]$ on $D$. Let $(\v_{s,t}\colon D\to D)$ be the  evolution family which solves the Loewner ODE $(\ref{solve})$.
If   $Q$ is a complex manifold of dimension $N$  and $(f_t\colon D\to Q)$ is a family of univalent mappings  satisfying  the functional equation $(\ref{functional})$ then $(f_t)$ solves the Loewner PDE $(\ref{LoewnerPDE})$. Furthermore there always exists a family of univalent mappings  $(f_t\colon D\to R)$ satisfying the functional equation $(\ref{functional})$ with values in an $N$-dimensional complex manifold $R$ which depends on $G(z,t)$.
\end{theorem}

The manifold $R$ is constructed as the direct limit of the evolution family $(\v_{s,t})$ associated with  $G(z,t)$ and satisfies $R=\cup_{t\geq 0}f_t(D)$. If $Q$ is an $N$-dimensional complex manifold and  $\Phi\colon R\to Q$ is a holomorphic mapping, then $(g_t:=\Phi\circ f_t)$ is a solution to the Loewner PDE (\ref{LoewnerPDE}), which is univalent if and only if $\Phi$ is univalent.

Conversely, if $(g_t:D\to Q)$ is a family of holomorphic mappings   which solves the  Loewner PDE  (\ref{LoewnerPDE}), then by the properties of the direct limit there exists a holomorphic map $\Phi: R\to Q$ such that $(g_t=\Phi\circ f_t)$, and the mapping $\Phi$ is univalent if and only if  $(g_t)$ is a univalent solution. Therefore if  $(g_t:D\to Q)$ is any univalent solution  to the Loewner PDE (\ref{LoewnerPDE}), the $N$-dimensional manifold $\cup_{t\geq 0}g_t(D)\subset Q$ is biholomorphic to $R$. For this reason the class of biholomorphism of $R$ is called  the {\sl Loewner range} of $G(z,t)$ (or of $(\v_{s,t})$).

The main open problem is thus the following: given a Herglotz vector field $G(z,t)$ of order $d\in[1,+\infty]$ on $D$, does there exist an univalent solution $(f_t:D\to \C^N)$  to the Loewner PDE (\ref{LoewnerPDE})? Or---given the previous discussion---equivalently, does the Loewner range $R$ of $G(z,t)$ embed as a domain of $\C^N$?

Note that, if $N=1$, then the Loewner range of any Herglotz vector field on the unit disc is an increasing union of discs, thus it is a simply connected open Riemann surface and by the uniformization theorem, it admits an embedding as a domain of $\C$.

In higher dimension there is no uniformization theorems available, and thus we need to study the union of increasing domains.

\section{The union problem}

Keeping in mind the discussion in Section \ref{chains}, we are going to study the embedding in $\C^N$  of $N$-dimensional complex manifolds $\Omega$  which are  the growing union of domains $\Omega_j$ biholomorphic to domains of $\C^N$. This problem turns out to arise naturally in complex dynamics, and it is related to the so called Bedford conjecture (see, {\sl e.g.} \cite{FornaessStensones} and \cite{Ar3}).

\begin{definition}
Let $D\subset D'\subset \C^N$ be two domains. The pair $(D,D')$ is a {\sl Runge pair} if $\mathcal{O}(D')$ is dense in $\mathcal{O}(D)$. A domain $D\subset \C^N$ is {\sl Runge} if $(D,\C^N)$ is a Runge pair.
\end{definition}

It is known \cite{ElKasimi} that  starlike  domains are Runge. Now we define some more general domains, which are easily seen to be Runge too.
\begin{definition}
We will say that a domain $D\subset\mathbb C^N$ is {\sl starshapelike}
if there exists an $\alpha\in\Aut_{hol}\mathbb C^N$ such that $\alpha(D)$
is starlike.
\end{definition}

 In the following, we will use this result:
\begin{theorem}[Anders\'{en}-Lempert, \cite{AndersenLempert}]\label{AL}
Let $D$ be a starshapelike domain in $\mathbb C^N$ for $N\geq 2$, and
let $\varphi: D\rightarrow\mathbb C^N$ be a biholomorphism onto a Runge
domain $\varphi(D)$.  Then there exists a sequence $\phi_j\in\Aut_{hol}\mathbb C^N$
such that $\phi_j\rightarrow\varphi$ uniformly on compact subsets of $D$.
\end{theorem}

The proof of our main theorem is based on the following result:

\begin{theorem}\label{union}
Let $\Omega=\cup_{j\in \N}\Omega_j$ be a complex manifold such that the following holds:
\begin{itemize}
\item[(1)] each pair $(\Omega_j,\Omega_{j+1})$ is a Runge pair, and
\item[(2)] each $\Omega_j$ is biholomorphic to a Stein starshapelike domain in $\mathbb C^N, N\geq 2$.
\end{itemize}
Then $\Omega$ is biholomorphic to a Runge and Stein domain in $\mathbb C^N$.
\end{theorem}
\begin{proof}
It is clear that $\Omega$ is Stein.  Choose Runge domains $U_j\subset\subset U_j'\subset\subset\Omega_j$ such that $\Omega=\cup_j U_j$,
and choose compact sets $K_j\subset\Omega_j$ with $U_j'\subset\subset K_j^\circ$.  Make
sure also that $U_j'\subset\subset U_{j+1}$.
Fix a sequence of biholomorphisms $(\varphi_j:\Omega_j\rightarrow\mathbb C^N)$ onto starshapelike domains.
We will inductively construct a sequence of embeddings $(\psi_j:\Omega_j\rightarrow\mathbb C^N)$ satisfying for all $j\in \N$ (for an $\epsilon$ to be specified),
\begin{itemize}
\item[$(0_j)$] $\psi_j(\Omega_j)$ is a starshapelike domain,
\item[($1_j$)] $\|\psi_{j+1} - \psi_{j}\|_{K_{j}}<\epsilon\cdot 2^{-j}$,
\item[($2_j$)] $\psi_j(U_{i-1})\subset\subset\psi_i(U_i')$ for $1\leq i\leq j$,
\item[$(3_j)$] $\psi_i(U_i')\subset\subset\psi_j(U_{i+1})$ for $0\leq i<j$.
\end{itemize}
We start by setting $\psi_0:=\varphi_0$.   Assume by inductive hypothesis that for all $0\leq j\leq m$
we have constructed embeddings $\psi_j:\Omega_j\rightarrow\mathbb C^N$ satisfying the previous conditions.

To construct $\psi_{m+1}$, notice that $\psi_m=\psi_m\circ\varphi_{m+1}^{-1}\circ\varphi_{m+1}$.
Consider now the biholomorphic mapping $\varphi_{m+1}\circ\psi_m^{-1}\colon \psi_m(\Omega_m)\to \varphi_{m+1}(\Omega_m)$. Its domain is starshapelike by assumption $0_m$, and its image is Runge since the pair $(\varphi_{m+1}(\Omega_m),\varphi_{m+1}(\Omega_{m+1}))$ is Runge by assumption $(1)$, and the domain $\varphi_{m+1}(\Omega_{m+1})$ is Runge, being starshapelike.
By  Theorem \ref{AL}   there exists a sequence $(\phi_{m+1}^s)\in\Aut_{hol}\mathbb C^N$ such that
$(\phi_{m+1}^s)^{-1}\rightarrow\varphi_{m+1}\circ\psi_m^{-1}$ uniformly
on $\psi_m(K_m)$ as $s\rightarrow\infty$.
Thus, it follows that $\phi_{m+1}^s\rightarrow\psi_m\circ\varphi_{m+1}^{-1}$ uniformly
on $\varphi_{m+1}(K_m)$ (see, {\sl e.g.}, \cite[Thm. 5.2]{DE}).

Therefore
$\psi_{m+1}^s:=\phi_{m+1}^s\circ\varphi_{m+1}\rightarrow\psi_m$ uniformly on
$K_m$ as $s\rightarrow\infty$, and so by choosing a large enough $\tilde{s}$ and defining $\psi_{m+1}:=\psi^{\tilde{s}}_{m+1}$, we get  $(1_{m+1})$.
It is easy to see that $(2_{m+1})$ and $(3_{m+1})$ also hold if $s$ is large enough and $(0_{m+1})$ is  true for any $s$.

Now it follows from $(1_m)$ that the sequence $\psi_j$ converges uniformly to
a map $\psi:\Omega\rightarrow\mathbb C^N$.  And if $\epsilon$ is chosen small enough,
depending on $\varphi_1$, it follows by Hurwitz that $\psi$ is injective.    By taking
the limit as $j\rightarrow\infty$ in $(2_m)$ and ($3_m$) we see that
$$
\psi(U_{i-1})\subset\psi_i(U_i')\subset\psi(U_{i+1}),
$$
for all $i\geq 1$, and so $\{\psi_i(U_i')\}$ is an exhaustion of $\psi(\Omega)$.  Since
each $\psi_{i}(U_i')$ is a Runge domain it follows that $\psi(\Omega)$ is Runge.
\end{proof}

\begin{remark}
A  result similar to Theorem \ref{union} was stated in \cite{Wold}, where each $\Omega_j$ is assumed to be biholomorphic to $\mathbb C^N$
and the conclusion is that $\Omega$ is biholomorphic to $\mathbb C^N$.
\end{remark}

\begin{remark}\label{refo}
Without  assumption (1), Theorem \ref{union} is false: indeed Forn\ae ss  \cite{Fornaess} gave an example of a non-Stein manifold $\Omega$ which is the growing union of domains $\Omega_j$ biholomorphic to the ball $\B^N$, $N=3$. In this example the pair $(\Omega_j,\Omega_{j+1})$ is not  Runge for all $j\in \N$.   In \cite{Wold} the third
author gave an example of a non-Stein manifold $\Omega$ which is the growing union domains $\Omega_j$ biholomorphic to $\mathbb C^2$.
\end{remark}

\begin{remark} In \cite{FornaessStensones}, Forn\ae ss and Stens\o nes proved that the direct limit  $\Omega$ (also called {\sl abstract basin of attraction} or {\sl tailspace}) of  a sequence
of univalent mappings $(\varphi_{j,j+1}:\mathbb B^N\rightarrow\mathbb B^N)_{j\in\N}$ all satisfying
\begin{equation}\label{conditionab}
a\|z\|\leq \|\varphi_{j,j+1}(z)\|\leq b\|z\|,
\end{equation}
for $0<a<b<1$, embeds as a Runge and Stein domain in $\mathbb C^N$.
It is not hard to see that condition (\ref{conditionab}) allows us to rescale the ball in such a way that $\varphi_{j,j+1}(\mathbb B^N)$ is  Runge for all $j\in \N$ (giving us an isomorphic direct limit). Denoting $(f_j\colon \B^N\to \Omega)$ the canonical morphisms, one has that $f_j=f_{j+1}\circ \v_{j,j+1}$ for all $j\in \N$, and thus $(f_j(\B^N),f_{j+1}(\B^N))$ is a Runge pair. Hence Forn\ae ss' and Stens\o nes' theorem  follows from Theorem \ref{union}.
\end{remark}

\section{The proof of Theorem \ref{main-intro}}
Recall the following definition by Docquier and Grauert \cite[Definition  20]{DocquierGrauert}.
\begin{definition}
Let $M$ be a nonempty open subset of a complex manifold $\tilde M$.
Then $M$ is {\sl semicontinuously holomorphically extendable to $\tilde M$ by means of a family $(M_t)_{0\leq t\leq 1}$ of nonempty open subsets of  $\tilde M$} iff the following holds:
\begin{itemize}
\item[(0)] $M_t$ is a Stein manifold for all $t$ in a dense subset of $[0,1]$,
\item[(1)] $M_0=M$ and $\bigcup_{0\leq t< 1}M_t=\tilde M$,
\item[(2)] $M_s\subset M_t$ for all $ 0\leq s\leq t\leq 1$,
\item[(3)] $\bigcup_{0\leq t< t_0}M_t$ is a union of connected components of $M_{t_0}$, for $0<t_0\leq 1$,
\item[(4)] $M_{t_0}$ is a union of connected components of the interior part of  $\bigcap_{t_0<t\leq 1}M_t$,  for $0<t_0\leq 1$.
\end{itemize}
\end{definition}

The following result is proved in \cite[Satz 17-Satz 19]{DocquierGrauert}.
\begin{theorem}\label{DG}
If $M$ and $\tilde M$ are Stein manifolds and $M$ is semicontinuously holomorphically extendable to $\tilde M$, then $(M,\tilde M)$ is a Runge pair.
\end{theorem}

Let $D\subset \C^N$ be a starlike complete hyperbolic domain. Let $G(z,t)$ be a Herglotz vector field of order $d\in[1,+\infty]$ on $D$.
By Theorem \ref{nostro} there exists a univalent solution $(f_t\colon D\to R)$ to the Loewner PDE (\ref{LoewnerPDE}), where $R=\cup_{t\geq 0}f_t(D)$ is the Loewner range of $G(z,t)$. According to what we stated in Section \ref{chains}, it is enough to prove that $R$ can be embedded as a Runge and Stein open subset of $\C^N$.

Recall that by Remark \ref{cresce} one has $f_s(D)\subset f_t(D)$ for all $0\leq s\leq t$.
We claim that for all $n\in \N$ the manifold $f_n(D)$ is semicontinuously holomorphically extendable to $f_{n+1}(D)$ by means of the family $(f_{n+t}(D))_{0\leq t\leq 1}$. Indeed (0),(1) and (2) are trivial, and (3) and (4) easily follow from i) of Definition \ref{solution} (note that the sequence of inverse maps
is a normal family due to hyperbolicity). Thus  Theorem \ref{DG} yields that $(f_n(D), f_{n+1}(D))$ is  a Runge  pair for any $n\in \N$.

Since $D$ is a Stein domain,  we can apply  Theorem \ref{union} to $R=\cup f_n(D)$ and we obtain that $R$ is biholomorphic to a Runge and Stein domain in $\mathbb C^N$.

\section{Embedding in evolution families}
Recall  \cite[Question 3]{AB}: if  $\v\colon \B^N\to\B^N$ is an univalent self-mapping of the unit ball, does there exists an
evolution family $(\v_{s,t}\colon \B^N\to\B^N)$ of order $d\in[1,+\infty]$ such that $\v_{0,1}=\varphi$?

The next proposition shows that without assuming that $\v(\B^N)$ is  Runge the answer is negative.
\begin{prop}
Let $D$ be a complete hyperbolic domain in $\mathbb{C}^N$.
Let  $(\v_{s,t}\colon D \to D)$ be an evolution family of order $d\in[1,+\infty]$. Then $(\v_{0,1}(D),D)$ is a Runge pair.
\end{prop}
\begin{proof}
Set $f_s:=\v_{s,1}$ for all $0\leq s\leq 1$.  Then $\v_{0,1}(D)$ is semicontinuously holomorphically extendable to $D$ by means of the family $(\v_{s,1}(D))_{0\leq s\leq 1}$.  Since $D$ is a Stein domain, Theorem \ref{DG} gives the result.
\end{proof}

This answers also a question in \cite[Section 9.4]{Ar1}:  we recall it briefly. With notations as in Remark \ref{refo},  let $f_j\colon \B^N\to \Omega_j$ be a biholomorphism and set $\v_{j,j+1}:=f_{j+1}^{-1}\circ f_j$ for all $j\in \N$.
Is the family $(\v_{j,j+1}\colon \B^N\to\B^N)$ embeddable in an evolution family $(\v_{s,t}\colon \B^N\to\B^N)$ of order $d\in[1,+\infty]$? The answer is negative: $\v_{0,1}(\B^N)$ is not Runge since by construction the pair $(\Omega_0,\Omega_1)$ is not Runge.

\bibliographystyle{amsplain}

\begin{thebibliography}{10}

\bibitem{ABCD} M. Abate, F. Bracci, M. D. Contreras, and S. D\'iaz-Madrigal, {\sl The evolution of Loewner's differential
equations}, Newsletter European Math. Soc. \textbf{78} (2010), 31--38.

\bibitem{AndersenLempert}
E. Anders\'{e}n,  and  L. Lempert, {\sl On the group of holomorphic automorphisms of $\mathbb C^n$},
Invent. Math.  \bf 110\rm (1992), no. 2, 371--388.

\bibitem{Ar1} L. Arosio, {\sl Resonances in Loewner equations},
 Adv.  Math. \textbf{227} (2011), 1413--1435.

\bibitem{Ar2}
L. Arosio, \textit{Loewner equations on complete hyperbolic domains}, preprint (arXiv:1102.5454 [math.CV])

\bibitem{Ar3} L. Arosio, {\sl Basins of attraction in Loewner equations}, Ann. Acad. Sci. Fenn. Math., to appear, doi:10.5186/aasfm.2012.3742, arXiv:1108.6000v2 [math.CV].


\bibitem{AB} L. Arosio and F. Bracci, {\sl Infinitesimal generators and the Loewner equation on complete hyperbolic manifolds}, Anal. Math. Phys., \textbf{1}  (2011), no. 4, 337--350.

\bibitem{ABHK} L. Arosio, F. Bracci, H. Hamada and G. Kohr, {\sl An abstract approach to Loewner chains}, J. Anal. Math., to appear.

\bibitem{BCD1} F. Bracci, M.D. Contreras and S.
D\'{\i}az-Madrigal, {\sl Evolution Families and the Loewner Equation I: the unit disc},  J. Reine Angew. Math. (Crelle's Journal), to appear, arXiv:0807.1594.

\bibitem{BCD2} F. Bracci, M.D. Contreras and S.
D\'{\i}az-Madrigal, {\sl Evolution Families and the Loewner Equation II: complex hyperbolic manifolds}, Math.
Ann. \textbf{344} (2009), 947--962.

\bibitem{one}
M. D. Contreras, S. D\'iaz-Madrigal and P. Gumenyuk, \textit{Loewner chains in the unit disc}, Rev. Mat. Iberoamericana, \textbf{26} (2010), 975--1012.

\bibitem{duality}
M. D. Contreras, S. D\'iaz-Madrigal and P. Gumenyuk, \textit{Local duality in Loewner equations}, preprint (arXiv:1202.2334v2 [math.DS]).

\bibitem{DE} P. G. Dixon, J. Esterle, \textit{Michael's problem and the Poincar\'e-Fatou-Bieberbach phenomenon}. Bull. Amer. Math. Soc. \textbf{15}, 2, (1986), 127---187.


\bibitem{DocquierGrauert}
F. Docquier and H. Grauert, {\sl Levisches Problem und Rungescher Satz f\"ur Teilgebiete Steinscher Mannigfaltigkeiten}.
Math. Ann. \bf 140\rm (1960), 94--123.

\bibitem{DGHK}
P. Duren, I. Graham, H. Hamada and G. Kohr \textit{Solutions for the generalized Loewner differential equation in several complex variables},
Math. Ann. \textbf{347} (2010), no. 2, 411--435.

\bibitem{ElKasimi}
A. El Kasimi, \textit{Approximation polyn\'{o}miale dans les domaines etoil\'es de $\C^n$}, Complex variables, \textbf{10} (1988), 179--182.

\bibitem{Fornaess}
J. E. Forn\ae ss, \textit{An Increasing Sequence of Stein Manifolds whose Limit is not Stein}, Math. Ann. \textbf{223} (1976), 275--277.

\bibitem{FornaessStensones}
J. E. Forn\ae ss and B. Stens\o nes,
{\sl Stable manifolds of holomorphic hyperbolic maps}, Internat. J. Math. \bf 15 \rm (2004), no. 8, 749--758.

\bibitem{GHK1}
I. Graham, H. Hamada and G. Kohr, {\sl Parametric representation
of univalent mappings in several complex variables}, Canadian J.
Math., \textbf{54} (2002), 324--351.

\bibitem{GHK2}I. Graham, H. Hamada, G. Kohr, and M. Kohr, {\sl Asymptotically
spirallike mappings in several complex variables}, J. Anal. Math.,
\textbf {105} (2008), 267--302.


\bibitem{Kuf1943} P.P. Kufarev, {\sl On one-parameter families of analytic functions,} (in Russian)
Mat. Sb. \textbf{13} (1943), 87--118.


\bibitem{Loewner} Ch. Loewner, {\sl Untersuchungen \"{u}ber schlichte
konforme Abbildungen des Einheitskreises}, Math. Ann.
\textbf{89} (1923), 103--121.

\bibitem{Pf74}
J.A. Pfaltzgraff, {\sl Subordination chains and univalence of
holomorphic mappings in $\mathbb{C}^n$}, Math. Ann.,
\textbf{210} (1974), 55--68.


\bibitem{Pf75}
J.A. Pfaltzgraff, {\sl Subordination chains and quasiconformal
extension of holomorphic maps in $\mathbb{C}^n$}, {Ann. Acad. Scie.
Fenn. Ser. A I Math.}, \textbf{1} (1975), 13--25.

\bibitem{Pommerenke-65}Ch. Pommerenke, {\sl \"{U}ber dis subordination
analytischer funktionen}, J. Reine Angew Math. \textbf{218} (1965), 159--173.


\bibitem{Por91}
T. Poreda, {\sl On generalized differential equations in Banach
spaces}, Dissertationes Mathematicae, \textbf{310} (1991), 1--50.

\bibitem{Voda}
M. Voda, \textit{Solution of a Loewner chain equation in several variables},  J. Math. Anal. Appl. \textbf{375} (2011), no. 1, 58--74.

\bibitem{Wold}
E. F. Wold, {\sl A long $\mathbb C^2$ which is not Stein},
Ark. Mat. \bf 48\rm  (2010), no. 1, 207--210.

\end{thebibliography}

\end{document}